\documentclass[12pt]{article}
\title{Realizations of automorphism groups of metric graphs induced by rational maps}
\author{Song JuAe \footnote{Tokyo Metropolitan University 1-1 Minami-Ohsawa, Hachioji, Tokyo, 192-0397, Japan. E-mail: song-juae@ed.tmu.ac.jp}}
\date{}

\usepackage{amsmath,amssymb,amscd}
\usepackage{amsthm}
\usepackage{color}
\usepackage{subfigure}
\usepackage{comment}

\newtheorem{dfn}{Definition}
\newtheorem{thm}[dfn]{Theorem}
\newtheorem{prop}[dfn]{Proposition}
\newtheorem{cor}[dfn]{Corollary}
\newtheorem{lemma}[dfn]{Lemma}
\newtheorem{rem}[dfn]{Remark}
\newtheorem{ex}[dfn]{Example}

\def\Gamma{\varGamma}

\begin{document}

\maketitle

\begin{abstract}
For a rational map $\phi$ from a metric graph $\Gamma$ to a tropical projective space $\boldsymbol{TP^n}$ defined by a ratio of rational functions $f_1, \ldots, f_{n + 1}$, an automorphism $\sigma$ of $\Gamma$ induces a permutation of the coordinates of $\boldsymbol{TP^n}$ if $\{ f_1, \ldots, f_{n + 1} \}$ is $\langle \sigma \rangle$-invariant.
Through this description, we can realize the automorphism group of $\Gamma$ as ambient automorphism group such as tropical projective general linear group, tropical general linear group and $\boldsymbol{Z}$-linear transformation group of Euclidean space.
\end{abstract}

{\bf keywords}: metric graphs, automorphism groups of metric graphs, rational maps, linear systems


\section{Introduction}

A {\it metric graph} $\Gamma$ is the underlying metric space of the pair of a graph $G$ and a length function $l: E(G) \rightarrow {\boldsymbol{R}}_{>0}$.
Here, a graph means an unweighted, undirected, finite, connected nontrivial multigraph and we allow the existence of loops and $E(G)$ denotes the set of edges of $G$.
In this paper, we give a way to realize (finite) automorphism groups of metric graphs as ambient automorphism groups such as tropical projective linear groups, tropical linear groups and $\boldsymbol{Z}$-linear transformation groups of Euclidean spaces.
We can simultaneously get these realizations by choosing one suitable set of rational functions, which is easy to find.
These realizations can be given by permutation matrices.

For any graph $G$, let $\boldsymbol{1}$ be the length function mapping all edge to one.
Then for the metric graph $\Gamma$ obtained from the pair $(G, \boldsymbol{1})$, we have a natural inclusion ${\rm Aut}(G) \hookrightarrow {\rm Aut}(\Gamma)$, where ${\rm Aut}(G)$ and ${\rm Aut}(\Gamma)$ denote the automorphism groups of $G$ and $\Gamma$, respectively.
${\rm Aut}(\Gamma)$ corresponds with the isometry transformation group of $\Gamma$ (cf. \cite{JuAe1} for the definition of (finite harmonic) morphism between metric graphs).
Thus, we can also have realizations of automorphism groups of graphs by our constructions.

Let $K$ be a complete algebraically closed non-Archimedean field with nontrivial valuation, and let $X, X^{\prime}$ be smooth, proper, connected curves over $K$.
If $\varphi : (X^{\prime}, V^{\prime} \cup D^{\prime}) \rightarrow (X, V \cup D)$ is a tame covering of triangulated punctured curve $(X, V \cup D)$ ($V, V^{\prime}$ are semistable vertex sets of $X, X^{\prime}$ respectively and $D \subset X(K)$, $D^{\prime} \subset X^{\prime}(K)$ punctures; see Definitions 3.8, 3.9, 4.25, 4.31 in \cite{ABBR1}) and the skeleton $\Sigma$ obtained from $(X, V \cup D)$ (see Subsection 3.7 of \cite{ABBR1}) has no loops, then the natural group homomorphism $\psi : {\rm Aut}_{X}(X^{\prime}) \rightarrow {\rm Aut}_{\Sigma}(\Sigma^{\prime})$ is injective by Theorem 7.4 (1) and Remark 7.5 in \cite{ABBR1}.
Here, $\Sigma^{\prime}$ denote the skeleton obtained from $(X^{\prime}, V^{\prime} \cup D^{\prime})$ and ${\rm Aut}_X(X^{\prime}), {\rm Aut}_{\Sigma}(\Sigma^{\prime})$ the automorphism groups of $\varphi$ and $\varphi|_{\Sigma}$, respectively.
(More precisely, see \cite{ABBR1} and \cite{ABBR2}.)
${\rm Aut}_X(X^{\prime})$ is a subset of the automorphism group of $(X^{\prime}, V^{\prime} \cup D^{\prime})$, so the automorphism group ${\rm Aut}(X^{\prime})$ of $X^{\prime}$.
And there is a natural group homomorphism $\psi^{\prime}$ from ${\rm Aut}_{\Sigma}(\Sigma^{\prime})$ to the automorphism group ${\rm Aut}(\Gamma^{\prime})$ of the underlying metric graph $\Gamma^{\prime}$ of $\Sigma^{\prime}$ (which may not be injective).
Therefore we can realize subgroups of ${\rm Aut}(X^{\prime})$ of the form ${\rm Aut}_X(X^{\prime})$ as the image of $\psi^{\prime} \circ \psi \subset {\rm Aut}(\Gamma^{\prime})$, and so as subgroup of our three groups.

We make an (injective) group homomorphism $\Psi$ from the automorphism group of a metric graph $\Gamma$ to the $\boldsymbol{Z}$-linear transformation group $\boldsymbol{Z}$-lin$(\boldsymbol{R^n})$ of $\boldsymbol{R^n}$ such that each automorphism of $\Gamma$ and the image by $\Psi$ are commutative with a rational map $\Gamma \rightarrow \boldsymbol{TP^n} \overset{i}\supset \boldsymbol{R^n}$, where $\boldsymbol{T}$ is the tropical semifield $(\boldsymbol{R} \cup \{ - \infty \}, {\rm max}, +)$, $\boldsymbol{TP^n}$ is the $n$-dimensional tropical projective space and $i$ is the inclusion $i : \boldsymbol{R^n} \hookrightarrow \boldsymbol{TP^n}; (X_1, \ldots, X_n) \mapsto (X_1 : \cdots : X_n : 0)$.
Note that we mean this inclusion $i$ whenever we write $\boldsymbol{R^n} \subset \boldsymbol{TP^n}$.
Concurrently, we also make other two realizations.
To make the group homomorphism, the following simple proposition is important.

\begin{prop}
	\label{main prop}
Let $\Gamma$ be a metric graph and $f_1, \ldots, f_{n + 1}$ distinct rational functions on $\Gamma$ other than the constant $- \infty$ function.
Let $\phi : \Gamma \rightarrow \boldsymbol{TP^n} ; x \mapsto (f_1(x) : \cdots : f_{n + 1}(x))$ be the induced rational map.
For $\sigma \in \rm{Aut}(\Gamma)$, if $\{ f_1, \ldots, f_{n + 1} \}$ is $\langle \sigma \rangle$-invariant, then $\sigma$ extends to a $\boldsymbol{Z}$-linear transformation of $\boldsymbol{R^n} \subset \boldsymbol{TP^n}$ through $\phi$, i.e., there is a regular $(n + 1) \times (n + 1)$ matrix $A_{\sigma}$ whose all coefficients are integers such that $i^{-1}(\phi(\sigma(x)))) = {}^t (A_{\sigma}({}^t i^{-1}(\phi(x))))$ holds for any $x \in \Gamma$.
\end{prop}

Here, ``$\{ f_1, \ldots, f_{n + 1} \}$ is $\langle \sigma \rangle$-invariant'' means that for any $k$, there exists a unique $l$ such that $f_k \circ \sigma = f_l$.
The following two propositions, other two cases we want, clearly hold since each permutation matrix is regular (see Subsection \ref{tropical algebra}).

\begin{prop}
	\label{main prop2}
Let $\Gamma$ be a metric graph and $f_1, \ldots, f_{n + 1}$ distinct rational functions on $\Gamma$ other than the constant $- \infty$ function.
Let $\phi : \Gamma \rightarrow \boldsymbol{TP^n} ; x \mapsto (f_1(x) : \cdots : f_{n + 1} (x))$ be the induced rational map.
For $\sigma \in {\rm Aut}(\Gamma)$, if $\{ f_1, \ldots, f_{n + 1} \}$ is $\langle \sigma \rangle$-invariant, then $\sigma$ extends to a permutation matrix in the tropical projective linear group ${\rm PGL}_{\rm trop} (n + 1, \boldsymbol{T})$.
\end{prop}

\begin{prop}
	\label{main prop3}
Let $\Gamma$ be a metric graph and $f_1, \ldots, f_{n + 1}$ distinct rational functions on $\Gamma$ other than the constant $- \infty$ function.
Let $\phi : \Gamma \rightarrow \boldsymbol{T^n} ; x \mapsto (f_1(x), \ldots, f_n (x))$ be the induced rational map.
For $\sigma \in {\rm Aut}(\Gamma)$, if $\{ f_1, \ldots, f_n \}$ is $\langle \sigma \rangle$-invariant, then $\sigma$ extends to a permutation matrix in the tropical linear group ${\rm GL}_{\rm trop} (n, \boldsymbol{T})$.
\end{prop}

$\boldsymbol{T^n}$ denotes the $n$-dimensional tropical affine space and see Subsection \ref{tropical algebra} for the definitions of ${\rm PGL}_{\rm trop} (n + 1, \boldsymbol{T})$ and ${\rm GL}_{\rm trop} (n, \boldsymbol{T})$.

By these propositions, our next goal is to find a way to get such rational functions.
As an answer, we use a complete linear system; from Proposition \ref{main prop}, we have the following corollary, which is the case that the rational map is induced by a complete linear system.

\begin{cor}
	\label{cor}
Let $\Gamma$ be a metric graph, $D$ a divisor on $\Gamma$.
For $\sigma \in {\rm Aut}(\Gamma)$, if the $\langle \sigma \rangle$-invariant linear system $|D|^{\langle \sigma \rangle}$ is not empty, then there exists a minimal generating set of $R(D)$ such that $\sigma$ extends to a $\boldsymbol{Z}$-linear transformation of $\boldsymbol{R^n} \subset \boldsymbol{TP^n}$ through the induced rational map $\Gamma \rightarrow \boldsymbol{TP^n}$.
\end{cor}

Here, $R(D)$ denotes the set of rational functions corresponding to elements of the complete linear system $|D|$ together with the constant $-\infty$ function and $|D|^{\langle \sigma \rangle}$ is $\{ D^{\prime} \in |D| \, | \, \forall x \in \Gamma, D^{\prime} (\sigma (x))= D^{\prime} (x) \}$ (which becomes a linear system; see Theorem 3.17 in \cite{JuAe1}).
From Propositions \ref{main prop2} and \ref{main prop3}, we have the following two corollaries:

\begin{cor}
	\label{cor2}
Let $\Gamma$ be a metric graph and $D$ a divisor on $\Gamma$.
For $\sigma \in {\rm Aut}(\Gamma)$, if the $\langle \sigma \rangle$-invariant linear system $|D|^{\langle \sigma \rangle}$ is not empty, there exists a minimal generating set of $R(D)$ such that $\sigma$ extends to a permutation matrix in ${\rm PGL}_{\rm trop}(n + 1, \boldsymbol{T})$ through the induced rational map $\Gamma \rightarrow \boldsymbol{TP^n}$.
\end{cor}

\begin{cor}
	\label{cor3}
Let $\Gamma$ be a metric graph and $D$ a divisor on $\Gamma$.
For $\sigma \in {\rm Aut}(\Gamma)$, if the $\langle \sigma \rangle$-invariant linear system $|D|^{\langle \sigma \rangle}$ contains an element $D^{\prime}$, there exists a $\langle \sigma \rangle$-invariant minimal generating set of $R(D^{\prime})$ such that $\sigma$ extends to a permutation matrix in ${\rm GL}_{\rm trop}(n, \boldsymbol{T})$ through the induced rational map $\Gamma \rightarrow \boldsymbol{T^n}$.\end{cor}

By these corollaries, we can realize subgroups of automorphism groups of metric graphs which is generated by one element.
Next, we prove each finite subgroup case; by Corollary \ref{cor}, we prove the following theorem:

\begin{thm}
	\label{main thm}
Let $\Gamma$ be a metric graph and $D$ a divisor on $\Gamma$.
Assume that the complete linear system $|D|$ induces an injective rational map $\Gamma \hookrightarrow \boldsymbol{TP^n}$.
For a finite subgroup $G$ of ${\rm Aut}(\Gamma)$, if the $G$-invariant linear system $|D|^G$ is not empty, then there exists a minimal generating set of $R(D)$ which induces an injective group homomorphism from $G$ to $\boldsymbol{Z}$-linear transformation group of $\boldsymbol{R^n} \subset \boldsymbol{TP^n}$ such that each element of $G$ and the image are commutative with the induced rational map (which may not be the original one).
\end{thm}

$|D|^G$ is the set $\{ D^{\prime} \in |D| \,|\, \forall x \in \Gamma, \forall \sigma \in G, D^{\prime}(\sigma (x)) = D^{\prime}(x) \}$ (and becomes a linear system by Theorem 3.17 in \cite{JuAe1} again).
Since for a metric graph $\Gamma$ which is not homeomorphic to a circle, ${\rm Aut}(\Gamma)$ is finite, by this theorem, we can realize it as a subgroup of $\boldsymbol{Z}$-lin$(\boldsymbol{R^n})$.
Other two cases are as follows:

\begin{thm}
	\label{main thm2}
Let $\Gamma$ be a metric graph and $D$ a divisor on $\Gamma$.
Assume that the complete linear system $|D|$ induces an injective rational map $\Gamma \hookrightarrow \boldsymbol{TP^n}$.
For a finite subgroup $G$ of ${\rm Aut}(\Gamma)$, if the $G$-invariant linear system $|D|^G$ is not empty, then there exists a minimal generating set of $R(D)$ which induces an injective group homomorphism from $G$ to ${\rm PGL}_{\rm trop}(n + 1, \boldsymbol{T})$ such that the image consists only of permutation matrices and each element of $G$ and the image are commutative with the induced rational map (which may not be the original one).
\end{thm}

\begin{thm}
	\label{main thm3}
Let $\Gamma$ be a metric graph and $D$ a divisor on $\Gamma$.
Assume that the complete linear system $|D|$ induces an injective rational map $\Gamma \hookrightarrow \boldsymbol{TP^n}$.
For a finite subgroup $G$ of ${\rm Aut}(\Gamma)$, if the $G$-invariant linear system $|D|^G$ contains an element $D^{\prime}$, then there exists a $G$-invariant minimal generating set of $R(D^{\prime})$ which induces an injective group homomorphism from $G$ to ${\rm GL}_{\rm trop}(n + 1, \boldsymbol{T})$ such that the image consists only of permutation matrices and each element of G and the image are commutative with the induced rational map $\Gamma \rightarrow \boldsymbol{T^{n + 1}}$.
\end{thm}

``A minimal generating set of $R(D)$ is $G$-invariant'' means that it is $\langle \sigma \rangle$-invariant for any $\sigma \in G$.
Since canonically $\boldsymbol{R^n} \subset \boldsymbol{T^n}$ and each $n \times n$ permutation matrix is in $\boldsymbol{Z}\text{-lin}(\boldsymbol{R^n})$, we have the following from Theorem \ref{main thm3}:

\begin{thm}
	\label{main thm4}
Let $\Gamma$ be a metric graph and $D$ a divisor on $\Gamma$.
Assume that the complete linear system $|D|$ induces an injective rational map $\Gamma \hookrightarrow \boldsymbol{TP^n}$.
For a finite subgroup $G$ of ${\rm Aut}(\Gamma)$, if the $G$-invariant linear system $|D|^G$ contains an element $D^{\prime}$, then there exists a $G$-invariant minimal generating set of $R(D^{\prime})$ which induces an injective group homomorphism from $G$ to $\boldsymbol{Z}\text{-lin}(\boldsymbol{R^{n + 1}})$ such that the image consists only of permutation matrices and each element of G and the image are commutative with the induced rational map $\Gamma \rightarrow \boldsymbol{R^{n + 1}}$.
\end{thm}

One advantage of Theorem \ref{main thm} compared to Theorem \ref{main thm4} is that the dimension of the Euclidean space in Theorem \ref{main thm} is that in Theorem \ref{main thm4} minus one.
One disadvantage of Theorem \ref{main thm} compared to Theorem \ref{main thm4} is that the image of group homomorphism in Theorem \ref{main thm4} consists only of permutation matrices but not in Theorem \ref{main thm}.

This paper is organized as follows.
Section 2 briefly reviews some basics of tropical algebra and of metric graphs including how to make rational maps induced by (complete) linear systems, which were given in \cite{Haase=Musiker=Yu}.
Proofs of Proposition \ref{main prop}, Corollaries \ref{cor}, \ref{cor2}, \ref{cor3}, Theorems \ref{main thm}, \ref{main thm2}, \ref{main thm3} are given in Section 3.
The section includes one corollary of Theorems \ref{main thm}, \ref{main thm2}, \ref{main thm3}, \ref{main thm4} and three examples of low genus metric graph cases.

{\bf Acknowledgements.} The author thanks my supervisor Masanori Kobayashi, Yuki Kageyama, Yasuhito Nakajima, Kohei Sato and Shohei Satake for helpful comments.
This work was supported by JSPS KAKENHI Grant Number 20J11910.

\section{Preliminaries}

In this section, we recall some basic facts of tropical algebra and of metric graphs and some results in \cite{Haase=Musiker=Yu} which we need later.

\subsection{Tropical algebra}
	\label{tropical algebra}

Let $\boldsymbol{T}$ be the algebraic system $(\boldsymbol{R} \cup \{ - \infty \}, {\rm max}, +)$.
We write the maximum operation ${\rm max}$ as $\oplus$, the ordinary addition $+$ as $\odot$, respectively.
$\boldsymbol{T}$ becomes a semifield with these two operations and is called {\it tropical semifield}.
As in the conventional algebra, we extend these two operations to matrices and vectors.
By the {\it $n + 1$ dimensional tropical (affine) space} $\boldsymbol{T^{n + 1}}$ and tropical scalar multiplication by $\boldsymbol{T^{\times}} = \boldsymbol{R}$, we can define the {\it $n$ dimensional tropical projective space} $\boldsymbol{TP^n}$ as $\boldsymbol{T^{n + 1}} / \boldsymbol{T^{\times}}$ as in the conventional algebra.
$I \in \boldsymbol{T^{n \times n}}$ denotes the identity matrix.
A tropical matrix $A \in \boldsymbol{T^{n \times n}}$ is {\it regular} or {\it invertible} if there exists a tropical matrix $B \in \boldsymbol{T^{n \times n}}$ such that $A \odot B = B \odot A = I$.
\cite{Cuninghame-Green} and \cite{Gaubert=Plus} show that the only tropical regular matrices are generalized permutation matrices.
Here, a {\it permutation matrix} is a matrix obtained by permuting the rows and/or the columns of the identity matrix and a {\it generalized permutation matrix} is the product of a diagonal matrix and a permutation matrix.
The {\it tropical general linear group} ${\rm GL}_{\rm trop}(n, \boldsymbol{T})$ is defined to be the set of all tropical regular square matrices of order $n$.
The {\it tropical projective linear group} ${\rm PGL}_{\rm trop}(n, \boldsymbol{T})$ is defined to be ${\rm GL}_{\rm trop}(n, \boldsymbol{T})$ modulo tropical scalar multiplication by $\boldsymbol{T^{\times}}$.

\subsection{Metric graphs and related basic facts}

Let $\Gamma$ be a metric graph.
The {\it genus} $g(\Gamma)$ of $\Gamma$ is its first Betti number.
We have the equality $g(\Gamma) = \# E(G) - \# V(G) + 1$, where $V(G)$ is the set of vertices of $G$, respectively, for any underlying graph $G$ of $\Gamma$.

Let ${\rm Div}(\Gamma)$ be the free abelian group generated by all points of $\Gamma$, i.e., ${\rm Div}(\Gamma) := \oplus_{x \in \Gamma} \boldsymbol{Z}x$.
An element of ${\rm Div}(\Gamma)$ is a {\it divisor} on $\Gamma$.
When $D$ is a divisor on $\Gamma$, the sum of all coefficients of $D$ is called the {\it degree} of $D$.
For a point $x \in \Gamma$, the degree of $D$ at $x$ is denoted by $D(x)$.
$D$ is {\it effective}, written by $D\ge0$, if all coefficients of $D$ are nonnegative.
The set of all points of $\Gamma$ where the coefficients of $D$ are not zero is called the {\it support} of $D$.

Let $f : \Gamma \rightarrow \boldsymbol{R} \cup \{-\infty\}$ be a continuous map.
$f$ is a {\it rational function} on $\Gamma$ if $f \equiv -\infty$ or $f$ is a piecewise $\boldsymbol{Z}$-affine function.
Let ${\rm Rat}(\Gamma)$ denote the set of all rational functions on $\Gamma$.
For $f, g \in {\rm Rat}(\Gamma)$ and $a \in \boldsymbol{T}$, we define {\it tropical sum} of $f$ and $g$, and {\it tropical scalar multiplication} of $f$ by $a$ as pointwise tropical operations, i.e., $(f \oplus g) (x) := {\rm max}\{ f(x), g(x) \}, (a \odot f) (x) := a + f(x)$ for any $x \in \Gamma$.
By these operations, ${\rm Rat}(\Gamma)$ becomes a tropical semimodule over $\boldsymbol{T}$.
Note that in fact we can define tropical multiplication on ${\rm Rat}(\Gamma)$ and this makes ${\rm Rat}(\Gamma)$ a tropical semiring over $\boldsymbol{T}$.
However, we need not this fact in this paper.

For $f \in {\rm Rat}(\Gamma)^{\times} = {\rm Rat}(\Gamma) \backslash \{ - \infty \}$ and $x \in \Gamma$, let ${\rm ord}_x(f)$ denote the sum of the outgoing slopes of $f$ at $x$.
The {\it principal divisor} ${\rm div}(f)$ defined by $f$ is $\sum_{x \in \Gamma}{\rm ord}_x(f) \cdot x$.
We define a relation $\sim$ on ${\rm Div}(\Gamma)$ as follows.
For $D_1, D_2 \in {\rm Div}(\Gamma)$, $D_1 \sim D_2$ if there exists $f \in {\rm Rat}(\Gamma)^{\times}$ such that ${\rm div}(f) = D_1 - D_2$.
This relation $\sim$ becomes an equivalence relation, which is called {\it linear equivalence}.
By the linear equivalence $\sim$, for a divisor $D$ on $\Gamma$, the {\it complete linear system} $|D|$ associated to $D$ is defined as the set of all effective divisors linearly equivalent to $D$.
Corresponding to the complete linear system $|D|$, we write $R(D)$ as the union $\{ f \in {\rm Rat}(\Gamma)^{\times} \,|\, D + {\rm div}(f) \ge 0 \} \cup \{ - \infty \}$.
Then $R(D)$ becomes a tropical subsemimodule over $\boldsymbol{T}$ of ${\rm Rat}(\Gamma)$ with the tropical sum and scalar multiplication (\cite[Lemma 4]{Haase=Musiker=Yu}).

It is not clear that $R(D)$ is finitely generated, however, in fact it is true.
In \cite{Haase=Musiker=Yu}, the authors proved that $R(D)$ is generated by the extremals and the set of all extremals is unique and finite up to the tropical scalar multiplication and a complete system of representatives is minimal (\cite[Corollary 9]{Haase=Musiker=Yu}).
Here, $f \in R(D)$ is called {\it extremal} if $g, h \in R(D), f = g \oplus h$ implies $f = g$ or $f = h$.
Extremals are characterized in the language of subgraphs:

\begin{lemma}[{\cite[Lemma 5]{Haase=Musiker=Yu}}]
	\label{key lemma}
A rational function $f$ is an extremal of $R(D)$ if and only if there are not two proper subgraphs $\Gamma_k$ (i.e. $\Gamma_k \not= \Gamma, \varnothing$) covering $\Gamma$ (i.e. $\Gamma_1 \cup \Gamma_2 = \Gamma$) such that each can fire on $D + {\rm div}(f)$.
\end{lemma}

Here, a {\it subgraph} of $\Gamma$ means a compact subset of $\Gamma$ with a finite number of connected components and a subgraph $\Gamma^{\prime}$ of $\Gamma$ {\it can fire} on a divisor $D$ if for any its boundary point $x_0$, the outdegree of $\Gamma^{\prime}$ at $x_0$ in $\Gamma$ is not greater than the coefficient of $D$ at $x_0$.
Note that by Lemma \ref{key lemma}, we can find all extremals of $R(D)$ when $\Gamma$ and $D$ are given concretely.
Especially, it suffices that we look into only subgraphs whose all boundary points are in the support of $D + {\rm div}(f)$ to check whether a rational function $f$ is an extremal of $R(D)$.

\begin{rem}
	\label{remark1}
\upshape{
Let $D \sim D^{\prime}$.
Then $n = m$ and $\{ D + {\rm div}(f_1), \ldots, D + {\rm div}(f_n)\} = \{ D^{\prime} + {\rm div}(g_1), \ldots, D^{\prime} + {\rm div}(g_n) \}$ hold for any minimal generating sets $\{ f_1, \ldots, f_n \}$ of $R(D)$ and $\{ g_1, \ldots, g_m \}$ of $R(D^{\prime})$.
In fact, as $D \sim D^{\prime}$, there exists a rational function $f \in {\rm Rat}(\Gamma)^{\times}$ such that $D^{\prime} = D + {\rm div}(f)$.
Therefore $R(D)$ is isomorphic to $R(D^{\prime})$ via $R(D) \rightarrow R(D^{\prime}); h \mapsto h - f$ (the inverse correspondence is given by $R(D^{\prime}) \rightarrow R(D); h \mapsto h + f$) and we have $D + {\rm div}(f_k) = D^{\prime} - {\rm div}(f) + {\rm div}(f_k) = D^{\prime} + {\rm div}(f_k - f)$ for any $k$.
Since $f_k$ is an extremal of $R(D)$, by Lemma \ref{key lemma}, there are not two proper subgraphs covering $\Gamma$ such that each can fire on $D + {\rm div}(f) = D^{\prime} + {\rm div}(f_k - f)$.
Thus $f_k - f$ is an extremal of $R(D^{\prime})$, and this means conclusions we wanted above.
}
\end{rem}

\begin{rem}
	\label{remark2}
\upshape{
For $\sigma \in {\rm Aut}(\Gamma)$, if $D$ is $\langle \sigma \rangle$-invariant (i.e. for any $x \in \Gamma$, $D(\sigma(x)) = D(x)$ holds), then an extremal $f$ of $R(D)$ is mapped by $\sigma$ to another extremal (possibly $f$ itself) of $R(D)$.
In fact, $f \circ \sigma$ is in $R(D)$ since $0 \le (D + {\rm div}(f))(\sigma (x)) = D (\sigma (x)) + ({\rm div}(f))(\sigma (x)) = D(x) + ({\rm div}(f \circ \sigma)) (x)$ hold for any $x \in \Gamma$.
If $f \circ \sigma$ is not an extremal of $R(D)$, then by Lemma \ref{key lemma}, there are two proper subgraphs $\Gamma_1$ and $\Gamma_2$ covering $\Gamma$ such that each can fire on $D + {\rm div}(f \circ \sigma)$.
The proper subgraphs $\sigma^{-1} (\Gamma_1)$ and $\sigma^{-1} (\Gamma_2)$ cover $\Gamma$ and each can fire on $D + {\rm div}(f)$, and this means that $f$ is not an extremal of $R(D)$ by Lemma \ref{key lemma} again.
}
\end{rem}

\begin{rem}
	\label{remark3}
\upshape{
For a finite subgroup $G$ of ${\rm Aut}(\Gamma)$, if $D$ is $G$-invariant (i.e. for any $\sigma \in G$, $x \in \Gamma, D(\sigma (x)) = D(x)$ holds), then there exists a $G$-invariant minimal generating set $\{ f_1, \ldots, f_n \}$ of $R(D)$.
In fact, it is enough to choose each $f_k$ as the maximum value is zero.
For any $\sigma \in G$, for any $k$, there exists a unique $l$ such that $f_k \circ \sigma = f_l$ since $f_k$ and $f_l$ have the same maximum value zero and by Remark \ref{remark2}.
Section 2 of \cite{Haase=Musiker=Yu} is also helpful to understand this argument.
Since every rational function is an (ordinary) sum of chip firing moves plus a constant by Lemma 2 of \cite{Haase=Musiker=Yu}, choosing the maximum value of $f_k$ as zero corresponds to choosing this constant as zero.
}
\end{rem}

For a divisor $D$ on $\Gamma$, there is a natural one-to-one correspondence between the complete linear system $|D|$ and the projection of $R(D)$, i.e., $\boldsymbol{P}R(D) = (R(D) \setminus \{-\infty\}) / \boldsymbol{T}^{\times}$.
Thus $|D|$ has a structure of finitely generated tropical projective space and induces a rational map from $\Gamma$ to a tropical projective space.
Concretely, for a minimal generating set $\{ f_1, \ldots, f_{n + 1} \}$ of $R(D)$, which all are extremals of $R(D)$, the {\it rational map} $\phi_{|D|} : \Gamma \rightarrow \boldsymbol{TP^n}$ {\it induced by} $|D|$ is given by the correspondence $x \mapsto (f_1(x) : \cdots : f_{n + 1}(x))$ for any $x \in \Gamma$.
Note that we use the ratio in tropical meaning and there is an arbitrariness of the choice of a minimal generating set $\{ f_1, \ldots, f_{n + 1} \}$ of $R(D)$.
Exchanging $\{ f_1, \ldots, f_{n + 1} \}$ to another minimal generating set of $R(D)$ induces a (classical) parallel translation of the image and a renumbering.
In other word, $|D|$ define a rational map up to the action of ${\rm PGL}_{\rm trop}(n + 1, \boldsymbol{T})$ on $\boldsymbol{TP^n}$.
We can always find a divisor whose complete linear system induces an injective rational map (cf. \cite[Theorem 45]{Haase=Musiker=Yu}).
We can define a distance function on the image of a rational map and with this distance function, an injective rational map {\it induced} by a complete linear system always becomes an isometry (see \cite{JuAe2}), but in this paper, we need not this fact.

\section{Main results}

In this section, we give proofs of our main results, their corollaries and some examples.

First, we give our proof of Proposition \ref{main prop}.

\begin{proof}[Proof of Proposition \ref{main prop}]
Since $\{ f_1, \ldots, f_{n + 1} \}$ is $\langle \sigma \rangle$-invariant, $\sigma$ induces a permutation of $\{ 1, \ldots, n + 1 \}$.
There is a number $s$ in $\{ 1, \ldots, n + 1\}$ such that $\sigma(s) = n + 1$.
Let $A_{\sigma} = (a_{k, l})_{ 1 \le k, l \le n}$ be the $n \times n$ matrix given by
\[
a_{k, l} :=
\begin{cases} 1 & \text{if $k \not= s$ and $l = \sigma(k)$,}\\
-1 & \text{if $l = \sigma(n + 1)$, and}\\
0 & \text{otherwise.}
\end{cases}
\]
Then, $\sigma$ and $A_{\sigma}$ are commutative with $\phi$.
\end{proof}

\begin{rem}
\upshape{
In the construction of the $n \times n$ matrix $A_{\sigma}$ in the above proof, we can see a peculiar phenomenon in the tropical world that we can make $A_{\sigma}$ as a $\boldsymbol{Z}$-linear transformation of $\boldsymbol{R^n}$ unlike classical case since the tropical division is the usual subtraction.
}
\end{rem}

We specify here that the proof of Proposition \ref{main prop} was inspired by that of Corollary 7.5 in \cite{Izhakian=Johnson=Kambites} and thank the authors for their great works.

\begin{proof}[Proof of Corollary \ref{cor}]
By the assumption, there is an element $D^{\prime} \in |D|^{\langle \sigma \rangle}$, and thus there is a rational function $f \in R(D) \setminus \{ - \infty \}$ such that $D^{\prime} = D + {\rm div}(f)$.
By Remark \ref{remark2} and tropical scalar multiplication, there exists a $\langle \sigma \rangle$-invariant minimal generating set $\{ g_1, \ldots, g_{n + 1} \}$ of $R(D^{\prime})$.
By Proposition \ref{main prop}, there is an $n \times n$ matrix $A_{\sigma}$ whose all coefficients are integers, and which and $\sigma$ are commutative with the induced rational map $\Gamma \rightarrow \boldsymbol{TP^n}; x \mapsto (g_1(x) : \cdots : g_{n + 1}(x))$.
For each $k$, let $f_k := g_k + f$.
By Remark \ref{remark1}, $\{ f_1, \ldots, f_{n + 1} \}$ is a minimal generating set of $R(D)$.
Since $(g_1 (x) : \cdots : g_{n + 1} (x)) = (f_1 (x) - f (x) : \cdots : f_{n + 1} (x) - f (x)) = (f_1 (x) : \cdots : f_{n + 1} (x))$ hold for any $x \in \Gamma$, $j(\phi(\sigma(x))) = {}^tA_{\sigma} {}^t(j(\phi(x)))$ holds with $\phi : \Gamma \rightarrow \boldsymbol{TP^n}; x \mapsto (f_1 (x) : \cdots : f_{n + 1}(x))$ and $j : {\rm Im}(i) \hookrightarrow \boldsymbol{R^n}; (X_1 : \cdots : X_{n + 1}) \mapsto (X_1 - X_{n + 1}, \ldots, X_n - X_{n + 1})$.
\end{proof}

Using $\{f_1, \ldots, f_{n + 1}\}$ (resp. $\{g_1, \ldots, g_{n + 1}\}$) in the above proof and Proposition \ref{main prop2} (resp. Proposition \ref{main prop3}), we have Corollary \ref{cor2} (resp. Corollary \ref{cor3}).

\begin{proof}[Proof of Theorem \ref{main thm}]
By Corollary \ref{cor}, Remark \ref{remark3} and the injectivity of induced rational map, we have the conclusion.
\end{proof}

\begin{proof}[Proof of Theorem \ref{main thm2}]
By Corollary \ref{cor2}, Remark \ref{remark3} and the injectivity of induced rational map, we have the conclusion.
\end{proof}

\begin{proof}[Proof of Theorem \ref{main thm3}]
By Remark \ref{remark3}, there exists a $G$-invariant generating set $\{ g_1, \ldots, g_{n + 1} \}$ of $R(D^{\prime})$.
Since $|D|$ induces an injective rational map $\Gamma \hookrightarrow \boldsymbol{TP^n}$, $\phi : \Gamma \rightarrow \boldsymbol{T^{n + 1}}; x \mapsto (g_1 (x), \ldots, g_{n + 1} (x))$ is also injective.
In fact, if $\phi (x) = \phi (y)$, then $g_k (x) = g_k (y)$ for all $k$, and so $(g_1 (x) : \cdots : g_{n + 1} (x)) = (g_1 (y) : \cdots : g_{n + 1} (y))$ holds.
Thus, we have $x = y$.
By Corollary \ref{cor3} and the injectivity, we have the conclusion.
\end{proof}

\begin{rem}
\upshape{
Except genus one leafless metric graph case, since in other cases all metric graphs have finite automorphism groups, we can always find a divisor satisfying the conditions of Theorems \ref{main thm}, \ref{main thm2}, \ref{main thm3} for $G = {\rm Aut}(\Gamma)$.
}
\end{rem}

\begin{rem}
\upshape{
In the proofs of Propositions \ref{main prop}, \ref{main prop2}, \ref{main prop3}, Corollaries \ref{cor}, \ref{cor2}, \ref{cor3} and Theorems \ref{main thm}, \ref{main thm2}, \ref{main thm3}, essentially we only use the $\langle \sigma \rangle$-invariance (or $G$-invariance) of the rational function set $\{ f_1, \ldots, f_{n + 1} \}$ defining rational map and the injectivity of rational map.
Thus, in this case, we use not a complete linear system but a linear subsystem. 
Moreover, actually we need not take $\{ f_1, \ldots, f_{n + 1} \}$ as a subset of a minimal generating set of $R(D)$.
However, this construction is very practical since a minimal generating set of $R(D)$ for a suitable divisor $D$ always has the two properties and we can easily find such $D$.   
}
\end{rem}

A metric graph is {\it hyperelliptic} if it has a divisor of degree two and of rank one.
Here, for a divisor $D$ on a metric graph, its {\it rank} is defined to be the minimum integer $s$ such that for some effective divisor $E$ of degree $s + 1$, the complete linear system associated to $D - E$ is empty.
The {\it canonical divisor} $K_{\Gamma}$ of a metric graph $\Gamma$ is the divisor on $\Gamma$ whose coefficient at each point $x$ is the valency of $x$ minus two, where the {\it valency} of $x$ is the number of connected components of $U \setminus \{ x \}$ for any sufficiently small connected neighborhood $U$ of $x$.

\begin{cor}
	\label{hyperelliptic cor}
Let $\Gamma$ be a metric graph of genus at least two.
If $\Gamma$ is not hyperelliptic, then the canonical linear system $|K_{\Gamma}|$ induces an injective rational map $\phi : \Gamma \rightarrow \boldsymbol{TP^n}$ and an injective group homomorphism $\Psi : {\rm Aut}(\Gamma) \hookrightarrow \boldsymbol{Z}${\rm -lin}$(\boldsymbol{R^n})$ such that $\phi$ commutes with each element $\sigma$ of ${\rm Aut}(\Gamma)$ and $\Psi(\sigma)$, where $n$ is the number of elements of a minimal generating set of $R(K_{\Gamma})$ minus one.
\end{cor}

\begin{proof}
As $\Gamma$ has genus at least two, ${\rm Aut}(\Gamma)$ is finite.
By Theorem 49 of \cite{Haase=Musiker=Yu}, the canonical map $\phi_{|K_{\Gamma}|}$ is injective.
Since $K_{\Gamma}$ is ${\rm Aut}(\Gamma)$-invariant, by Theorem \ref{main thm}, we get our conclusion.
\end{proof}

By the same proof of Corollary \ref{hyperelliptic cor} using Theorems \ref{main thm2}, \ref{main thm3}, \ref{main thm4} instead of Theorem \ref{main thm}, respectively, we have the following three corollaries:

\begin{cor}
	\label{hyperelliptic cor2}
Let $\Gamma$ be a metric graph of genus at least two.
If $\Gamma$ is not hyperelliptic, then the canonical linear system $|K_{\Gamma}|$ induces an injective rational map $\phi : \Gamma \hookrightarrow \boldsymbol{TP^n}$ and an injective group homomorphism $\Psi : {\rm Aut}(\Gamma) \hookrightarrow {\rm PGL}_{\rm trop}(n, \boldsymbol{T})$ such that $\phi$ commutes with each element $\sigma$ of ${\rm Aut}(\Gamma)$ and $\Psi(\sigma)$, where $n$ is the number of elements of a minimal generating set of $R(K_{\Gamma})$ minus one.
\end{cor}

\begin{cor}
	\label{hyperelliptic cor3}
Let $\Gamma$ be a metric graph of genus at least two.
If $\Gamma$ is not hyperelliptic, then the canonical linear system $|K_{\Gamma}|$ induces an injective rational map $\phi : \Gamma \hookrightarrow \boldsymbol{T^n}$ and an injective homomorphism $\Psi : {\rm Aut}(\Gamma) \hookrightarrow {\rm GL}_{trop}(n, \boldsymbol{T})$ such that $\phi$ commutes with each element $\sigma$ of ${\rm Aut}(\Gamma)$ and $\Psi(\sigma)$, where $n$ is the number of elements of a minimal generating set of $R(K_{\Gamma})$.
\end{cor}

\begin{cor}
	\label{hyperelliptic cor4}
Let $\Gamma$ be a metric graph of genus at least two.
If $\Gamma$ is not hyperelliptic, then the canonical linear system $|K_{\Gamma}|$ induces an injective rational map $\phi : \Gamma \hookrightarrow \boldsymbol{R^n}$ and an injective homomorphism $\Psi : {\rm Aut}(\Gamma) \hookrightarrow \boldsymbol{Z}\text{-lin}(\boldsymbol{R^n})$ such that $\phi$ commutes with each element $\sigma$ of ${\rm Aut}(\Gamma)$ and $\Psi(\sigma)$, where $n$ is the number of elements of a minimal generating set of $R(K_{\Gamma})$.
\end{cor}

\begin{rem}
\upshape{
We can also make the projective space $\boldsymbol{RP^n}$ (in the usual sense) versions in the same arguments as up untile now since $\boldsymbol{R^{n + 1}} \subset \boldsymbol{RP^n}$ and by the definition of projective linear group.
Moreover, for a topological space $X$ (plus some additional structures) with its automorphism group ${\rm Aut}(X)$ (for a definition of automorphism of $X$), if $X$ contains $\boldsymbol{R^n}$ and ${\rm Aut}(X)$ contains all permutation matrices or elements of the form $A_{\sigma}$ in the proof of Proposition \ref{main prop} (or corresponding automorphisms), then we have the same conclusions for $X$.
}
\end{rem}

\begin{ex}
	\label{example1}
\upshape{
Let $\Gamma$ be the closed interval $[0, 1]$.
We call the point $0$ (resp. $1$, $1/2$) as $x$ (resp. $y$, $z$).
Let $\iota$ be the unique nontrivial automorphism of $\Gamma$, i.e., $\iota$ is an isometry $\Gamma \rightarrow \Gamma$ such that $\iota(x) = y$ holds.
We have ${\rm Aut}(\Gamma) = \langle \iota \rangle$.
Let $D = x$.
Then the rational function $f_1$ with slope one on $\Gamma$ and setting $f_1 (x) := 1$ and $f_1(y) := 0$ and the constant zero function $f_2$ on $\Gamma$ generate $R(D)$.
The image of $\Gamma \rightarrow \boldsymbol{TP^1} \supset \boldsymbol{R^1} ; x \mapsto (f_1 (x) : f_2 (x))$ is the closed interval $[0, 1] \subset \boldsymbol{R^1}$ and so $\iota$ induces a $\boldsymbol{Z}$-affine transformation of $\boldsymbol{R^1}$ but not a $\boldsymbol{Z}$-linear transformation of $\boldsymbol{R^1}$.
Since $\iota$ fixes $z$ and the $\langle \iota \rangle$-invariant linear system $|D|^{\langle \iota \rangle}$ contains the divisor $z=: D^{\prime}$, we can find an $\langle \iota \rangle$-invariant generating set $\{ f_1^{\prime}, f_2^{\prime} \}$ of $R(D^{\prime})$ such that $f_1^{\prime}|_{[ x, z]} \equiv 0$, $f_1^{\prime} (y) := -1/2$, $f_1^{\prime}$ has slope one on $[z, y]$ and $\iota (f_1^{\prime}) = f_2^{\prime}$ holds.
Then the induced rational map $\phi^{\prime} := (f_1^{\prime} : f_2^{\prime}) : \Gamma \rightarrow \boldsymbol{TP^1} \supset \boldsymbol{R^1}$ has the image $[ -1/2, 1/2] \subset \boldsymbol{R^1}$ and $\iota$ induces the square matrix $A_{\iota} = ( -1 )$.
Finally, we have the injective group homomorphism ${\rm Aut}(\Gamma) \hookrightarrow \boldsymbol{Z}${\rm -lin}$(\boldsymbol{R^1}); {\rm id}_{\Gamma} \mapsto (1), \iota \mapsto (-1)$, where ${\rm id}_{\Gamma}$ denotes the identity map of $\Gamma$.
Also $\phi^{\prime}$ induces ${\rm Aut}(\Gamma) \hookrightarrow {\rm PGL}_{\rm trop}(2, \boldsymbol{T}) ; {\rm id}_{\Gamma} \mapsto \left( \begin{array}{cc} - \infty & 0 \\ 0 & -\infty \end{array}\right), \iota \mapsto \left( \begin{array}{cc} - \infty & 0 \\ 0 & - \infty \end{array}\right)$ and $\Gamma \hookrightarrow \boldsymbol{T^2} ; x \mapsto (f_1(x), f_2(x))$ induces ${\rm Aut}(\Gamma) \hookrightarrow {\rm GL}_{\rm trop}(2, \boldsymbol{T}) ; {\rm id}_{\Gamma} \mapsto \left( \begin{array}{cc} - \infty & 0 \\ 0 & -\infty \end{array}\right), \iota \mapsto \left( \begin{array}{cc} - \infty & 0 \\ 0 & - \infty \end{array}\right)$.
}
\end{ex}

\begin{ex}
	\label{example2}
\upshape{
In the same setting as Example \ref{example1}, let $E:= 2z$.
Then $R(E)$ is generated by the three rational functions $g_1, g_2$ and $g_3$, where $g_1 |_{[x, z]} \equiv 0$, $g_1 (y) := -1$, $g_1$ has slope two on $[z, y]$, and $\iota (f_1) = f_2$ holds, and $g_3(x) := g_3(y) := -1/2$, $g_3 (z) := 0$ and $g_3$ has slope one on $[x, z]$ and $[z, y]$.
Since the set $\{ g_1, g_2, g_3 \}$ is $\langle \iota \rangle$-invariant and the induced rational map $\psi := (g_1 : g_2 : g_3) : \Gamma \rightarrow \boldsymbol{TP^2} \supset \boldsymbol{R^2}$ is injective, we have the injective group homomorphism ${\rm Aut}(\Gamma) \hookrightarrow \boldsymbol{Z}${\rm -lin}$(\boldsymbol{R^2}) ; {\rm id}_{\Gamma} \mapsto \left( \begin{array}{cc} 1 & 0 \\ 0 & 1 \end{array}\right), \iota \mapsto \left( \begin{array}{cc} 0 & 1 \\ 1 & 0 \end{array}\right)$.
Also $\psi$ induces ${\rm Aut}(\Gamma) \hookrightarrow {\rm PGL}_{\rm trop}(3, \boldsymbol{T}) ; {\rm id}_{\Gamma} \mapsto \left( \begin{array}{ccc} 0 & - \infty & - \infty \\ - \infty & 0 & -\infty \\ -\infty & -\infty & 0 \end{array}\right), \iota \mapsto \left( \begin{array}{ccc} - \infty & 0 & -\infty \\ 0 & - \infty & - \infty \\ - \infty & - \infty & 0 \end{array}\right)$ and $\Gamma \hookrightarrow \boldsymbol{T^3} ; x \mapsto (g_1(x),  g_2(x), g_3(x))$ induces ${\rm Aut}(\Gamma) \hookrightarrow {\rm GL}_{\rm trop}(3, \boldsymbol{T}) ; {\rm id}_{\Gamma} \mapsto \left( \begin{array}{ccc} 0 & - \infty & - \infty \\ - \infty & 0 & -\infty \\ - \infty & - \infty & 0 \end{array}\right), \iota \mapsto \left( \begin{array}{ccc} - \infty & 0 & - \infty \\ 0 & - \infty & -\infty \\ - \infty & - \infty & 0 \end{array}\right)$.
}
\end{ex}

\begin{ex}
	\label{example3}
\upshape{
Let $\Gamma$ be a circle of length four.
Fix a point $x \in \Gamma$.
Let $\sigma$ be the $180$ degrees rotation and $x^{\prime} := \sigma(x)$.
For the divisor $D := x + x^{\prime}$, we can choose a $\langle \sigma \rangle$-invariant minimal generating set $\{ f_1, f_2 \}$ of $R(D)$.
Concretely, if we call the midpoints of the two pathes $P_1$ and $P_2$ between $x$ and $x^{\prime}$ as $p_1$ and $p_2$ respectively, then for example, we can choose $f_1$ as $f_1 (x) := f_2 (x) := 1$, $f_1 (p_1) := 0$, $f_1$ has slope one on $[x, p_1] \cap P_1$ and $[x^{\prime}, p_1] \cap P_1$ and $f_1 |_{P_1} \equiv 1$, and $f_2 := f_1 \circ \sigma$.
Then, $\phi : \Gamma \rightarrow \boldsymbol{TP^1} \supset \boldsymbol{R^1}$ is not injective and the image in $\boldsymbol{R^1}$ is $[-1, 1]$.
$\phi$ induces the injective group homomorphisms $\langle \sigma \rangle \hookrightarrow \boldsymbol{Z}\text{-lin} (\boldsymbol{R^1}) ; {\rm id}_{\Gamma} \mapsto (1), \sigma \mapsto (-1)$, $\langle \sigma \rangle \hookrightarrow {\rm PGL}_{\rm trop}(2, \boldsymbol{T}); {\rm id}_{\Gamma} \mapsto \left( \begin{array}{cc} 0 & - \infty \\ -\infty & 0 \end{array}\right), \sigma \mapsto \left( \begin{array}{cc} - \infty & 0 \\ 0 & -\infty \end{array}\right)$.
Also $\Gamma \rightarrow \boldsymbol{T^2}; x \mapsto (f_1 (x), f_2 (x))$ induces the injective group homomorphism $\langle \sigma \rangle \hookrightarrow {\rm GL}_{\rm trop}(2, \boldsymbol{T}); {\rm id}_{\Gamma} \mapsto \left( \begin{array}{cc} 0 & - \infty \\ -\infty & 0 \end{array}\right), \sigma \mapsto \left( \begin{array}{cc} - \infty & 0 \\ 0 & -\infty \end{array}\right)$.
On the other hand, for the isometry $\iota : \Gamma \rightarrow \Gamma$ which maps each point to the line symmetric point with the line $xx^{\prime}$ as the axis of symmetry, $\{ f_1, f_2 \}$ is also $\langle \iota \rangle$-invariant, so we have three injective group homomorphisms from $\langle \iota \rangle$ having the same images as above.
Thus, $\phi$ and $(f_1, f_2)$ do not induce injective group homomorphisms from $\langle \sigma, \iota \rangle$ to $\boldsymbol{Z}\text{-lin}(\boldsymbol{R^1})$ or ${\rm PGL}_{\rm trop}(2, \boldsymbol{T})$ or ${\rm GL}_{\rm trop}(2, \boldsymbol{T})$.
It comes from the fact that $|D|$ dose not induce an injective rational map.
}
\end{ex}

\end{document}